\documentclass{amsart}
\usepackage{amsfonts,amssymb,amsmath,amsthm}
\usepackage{url}
\usepackage{enumerate}
\usepackage{bbm}
\usepackage{amssymb}

\newtheorem{thrm}{Theorem}[section]
\newtheorem{lem}[thrm]{Lemma}
\newtheorem{prop}[thrm]{Proposition}
\newtheorem{cor}[thrm]{Corollary}
\theoremstyle{definition}

\newtheorem{remark}[thrm]{Remark}
\numberwithin{equation}{section}

\DeclareMathOperator{\conv}{conv}

\newcommand{\R}{\ensuremath{{\mathbb R}}}
\newcommand{\E}{\ensuremath{{\mathbb E}}}
\newcommand{\norm}[1]{\left \lVert#1 \right\rVert}
\newcommand{\abs}[1]{\left\lvert#1 \right\rvert}
\newcommand{\skp}[1]{\left<#1\right>}

\author{David Alonso-Guti\'errez  \and Joscha Prochno}
\address[David Alonso-Guti\'errez]{Department of Mathematical and Statistical Sciences\\
University of Alberta\\
505 Central Academic Building\\
Edmonton T6G 2G1\\
Canada}
\email{alonsogu@ualberta.ca}
\address[Joscha Prochno]{Department of Mathematical and Statistical Sciences\\
University of Alberta\\
605 Central Academic Building\\
Edmonton T6G 2G1\\
Canada}
\email{prochno@ualberta.ca}
\thanks{~}

\begin{document}

\title[Gaussian Behavior and Mean width of random polytopes]{On the Gaussian behavior of marginals and the mean width of random polytopes}

\begin{abstract}
We show that the expected value of the mean width of a random
polytope generated by $N$ random vectors ($n\leq N\leq e^{\sqrt n}$) uniformly distributed
in an isotropic convex body in $\R^n$ is of the order $\sqrt{\log N} L_K$. This
completes a result of  Dafnis, Giannopoulos and Tsolomitis. 
We also prove some results in connection with the 1-dimensional marginals of the uniform probability measure on an isotropic convex body, extending the interval  
in which the average of the distribution functions
of those marginals behaves in a sub-
or supergaussian way.
\end{abstract}
\maketitle

\section{Introduction} \label{sect1}

In asymptotic convex geometry, the hyperplane conjecture is a very well known problem
that first appeared explicitly in \cite{key-Bo}. This conjecture says that there exists an absolute
constant $c$ such that every convex body of volume 1 has a hyperplane
section of volume greater than $c$. A result by Hensley \cite{key-H} yields an equivalent formulation, saying that there
exists an absolute constant $c$ such that every convex body has isotropic constant $L_K$ less than $c$.

The study of random polytopes began with Sylvester and the famous
four-point problem nearly 150 years ago. Since then, a tremendous
effort has been made to study expectations, variances, and
distributions of several functionals on a random polytope. This turned out to be very useful and many applications have been found (see \cite{key-Ba}, \cite{key-Re} and references therein).  Random polytopes also
provided counterexamples to several conjectures (see, for instance, \cite{key-G}, \cite{key-Sza} or \cite{key-N}).

In 1989, Milman and Pajor \cite{key-MP} showed a deep connection between the hyperplane conjecture
and the study of random polytopes, proving that the expected volume of
a random simplex in an isotropic convex body is closely related to the
value of its isotropic constant.

In \cite{key-DGT}, the authors studied the expected value of the
querma\ss integrals of a random polytope generated by $N$ random
vectors uniformly distributed in an isotropic convex body in $\R^n$. They showed that if
$n\leq N\leq e^{\sqrt n}$ the expected value of the smallest querma\ss integral, which is the volume
radius, is greater than $c\sqrt{\log\frac{N}{n}}L_K$ and the expected value of the biggest
one, which is the mean width, is smaller than $C\sqrt{\log
  N}L_K$. This yields a sharp estimate for the expected value of any
querma\ss integral when $n^2\leq N\leq e^{\sqrt n}$, but leaves a gap for the
range $n\leq N\leq n^2$. Our first purpose is to fill this gap for the
expected value of the mean width. Denoting by $a\sim b$ the fact that
there exist positive absolute constants $c, C$ such that $ca\leq b\leq
Ca$, in Section \ref{SEC expected value of mean width} we will prove the following:

\begin{thrm}\label{THM mean width}
Let $K\subseteq\R^n$ be a symmetric isotropic convex body and
$n\leq N\leq e^{\sqrt n }$. Let $K_N=\conv \{\pm X_1,\dots,\pm X_N\}$ be
a random polytope, where $X_1,\dots,X_N$ are independent random vectors
uniformly distributed in $K$. Then,
$$
\E w(K_N)\sim\sqrt{\log N}L_K.
$$
\end{thrm}

This estimate in the range $n\leq N\leq n^2$ will be a consequence of
the central limit theorem for convex bodies proved by Klartag \cite{key-K3} and the results proved by Sodin \cite{key-S}. The central limit theorem for convex
bodies was first considered in \cite{key-ABP} and says that most of the 1-dimensional marginals $\langle X,\theta\rangle$ of a random
variable $X$ uniformly distributed in an isotropic convex body $K$ are, in a certain sense, approximately Gaussian. To be more precise, Klartag showed that the distribution function $F_\theta(t)$ of most of
these marginals is ``almost'' Gaussian whenever $|t|$ is smaller than some power of
$n$. It turns out that the Gaussian behavior for a particular value of $t$ in this range will be enough to prove
Theorem \ref{THM mean width}.

Along these lines a great deal of research was devoted in connection with the marginals of the uniform probability measure on an isotropic convex body and important results were obtained. For
instance, in \cite{key-Bo2} Bourgain verified the hyperplane conjecture for the class of
$\psi_2$ bodies, {\it i.e.}, the class of convex bodies such that
every direction (or 1-dimensional marginal) is subgaussian.
Given an isotropic convex body $K$, we say that a direction $\theta\in S^{n-1}$ is
subgaussian with constant $r>0$ if
$$
|\{x\in K : |\langle x,\theta\rangle|\geq tL_K\}| \leq e^{-\frac{t^2}{r^2}}
$$
for all $1\leq t\leq r\sqrt n$.

In this setting, the following question was posed by Milman: is it
true that every convex body has at least one subgaussian direction?
This question has been answered in the affirmative for the class of
1-unconditional convex bodies \cite{key-BN2}, for the class of zonoids
\cite{key-P3}, and for the class of isotropic convex bodies with small diameter
\cite{key-P4}. (In fact, it was shown that the measure of subgaussian directions is greater than $1-e^{-\sqrt{n}}$ for this last class of convex bodies) In \cite{key-K2}, Klartag established the existence of a subgaussian
direction up to a logarithmic factor in the dimension. See also
\cite{key-GPV}.

In \cite{key-Pi}, Pivovarov considered the dual question of finding
supergaussian directions. We say that a direction $\theta\in S^{n-1}$
is supergaussian with constant $r>0$ if for all $1\leq
t\leq\frac{\sqrt n}{r}$ we have
$$
|\{x\in K\,:\,|\langle x,\theta\rangle|\geq tL_K\}|\geq e^{-r^2t^2}.
$$
He gave an affirmative answer up to a logarithmic factor for the class of 1-unconditional convex bodies.

In \cite{key-P5}, Paouris showed that every isotropic convex body with
bounded isotropic constant has
``many'' supergaussian directions. This includes several classes of
convex bodies such as 1-unconditional convex bodies, zonoids, duals of
zonoids, and the unit balls of the Schatten classes. He also proved that if for every
isotropic convex body a random direction is supergaussian with high
probability, then the hyperplane conjecture is true.

Going in the same direction, Klartag proved in \cite{key-K} that every
non-degenerate $n$-dimensional measure has one direction that behaves
in a ``supergaussian way'' for $t$ in the interval $1\leq t\leq c(\log
n)^\frac{1}{4}$. As a consequence of the aforementioned central limit
theorem this interval was extended to $1\leq t\leq n^\kappa$ for some
constant $\kappa$ (see \cite{key-S}).

The approach used to prove the central limit theorem, as well as some
previous weaker results (see \cite{key-ABP}, \cite{key-Bob} or \cite{key-S}), involved the study of the average of the
distribution functions of the 1-dimensional marginals together with
a concentration of measure phenomenon. More precisely,
Sodin proved that if the Euclidean norm verifies a
concentration hypothesis, then the average of the distribution
function of the 1-dimensional marginals of an isotropic random vector
is approximately Gaussian for $|t|\leq n^\kappa$ for some absolute
constant $\kappa$, which is smaller than $\frac{1}{4}$ (due to ``spherical approximation"). A similar result was obtained for $k$-dimensional
marginals in \cite{key-BB}. In \cite{key-K3}, the concentration hypothesis was shown to be
true for every convex body.

In \cite{key-AGP}, a new approach to study the expected value of the support
function of a random polytope generated by $N$ vertices in an
isotropic convex body was introduced.
Using this approach, we will prove in Section \ref{SuperAndSubGaussianAverage} that a supergaussian
estimate for the average of the distribution function of the 1-dimensional
marginals of the uniform probability measure on an isotropic body holds for the whole
range $1\leq t\leq\frac{\sqrt n}{c}$, and not only if $t\leq
n^\kappa$. Namely, if we define
$$
F(t)=\int_{S^{n-1}}|\{x\in K\,:\,|\langle x,\theta\rangle|\geq tL_K\}| \, d\sigma(\theta),
$$
where $d\sigma$ denotes the uniform probability measure on $S^{n-1}$,
we have the following:

\begin{thrm}\label{SupergaussianAverage}
There exists an absolute constant $c$ such that for every symmetric isotropic
convex body $K\subseteq\R^n$ and every $1\leq t\leq\frac{\sqrt n}{c}$ we have
$$
F(t)\geq e^{-c^2t^2}.
$$
\end{thrm}

Using the same idea, in Section \ref{SubgaussianEstimates} we will also show a subgaussian estimate for
the average distribution function. However, in this case the estimate
does not cover the whole range $1\leq t\leq c\sqrt n$ except for
convex bodies with small diameter, {\it i.e.}, $R(K)\leq C\sqrt n L_K$.

\begin{thrm}\label{SubgaussianAverage}
There exists an absolute constant $c$ such that for every symmetric isotropic
convex body $K\subseteq\R^n$ and any $1\leq t\leq n^\frac{1}{4}$
$$
F(t)\leq e^{-c^2t^2}.
$$
Furthermore, if $K$ has small diameter then this estimate is true for $1\leq t\leq
C\sqrt n$.
\end{thrm}

As a consequence, we will find an interval in which a random direction
$\theta\in S^{n-1}$ verifies a subgaussian estimate with high probability. A
restriction of the interval leads to better estimates for the measure of set of
``subgaussian directions'' in the case of convex bodies with small diameter.

\section{Preliminaries and Notation}

Before we go into more detail we start with some basic definitions. A convex body $K\subset \R^n$ is a compact convex set with non-empty interior. It is called symmetric if $-x\in K$, whenever $x\in K$. 
We will denote its volume (or Lebesgue measure) by $|\cdot |$. The volume of the Euclidean unit ball $B_2^n$ will be denoted by $w_n=|B_2^n|$. We write $S^{n-1} = \{ x \in \R^n : \norm{x}_2 = 1\}$ for the standard Euclidean sphere in $\R^n$ and $d\sigma$ for the uniform probability measure on $S^{n-1}$. A convex body is said to be isotropic if it has volume $1$ and satisfies the following two conditions:
\begin{itemize}
\item$\int_Kx \, dx=0 \textrm{ (center of mass at 0)}$,
\item$\int_K\langle x,\theta\rangle^2 \, dx=L_K^2\quad \forall \theta\in S^{n-1}$,
\end{itemize}
where $L_K$ is a constant independent of $\theta$, which is called the isotropic constant of $K$. Here, $\skp{\cdot,\cdot}$ denotes the standard scalar product in $\R^n$. 

A probability measure $\mu$ on $\R^n$ is said to be isotropic if it is
centered at $0$ and its covariance matrix is the identity. Notice
that a convex body is isotropic if and only if the uniform probability
measure on $\frac{K}{L_K}$ is isotropic.

Let $K$ be a convex body and $\theta\in S^{n-1}$ a unit vector. The support function of $K$ in the direction $\theta$ is defined by $h_K(\theta)=\max\{\skp{x,\theta}: x\in K\}$. The mean width of $K$ is 
$$
w(K)=\int_{S^{n-1}}h_K(\theta) \, d\sigma(\theta).
$$

Given a symmetric isotropic convex body $K$, we denote by $K_N=\textrm{conv}\{\pm X_1,\dots,\pm X_N\}$ the random polytope, where $X_1,\dots, X_N$ are independent random vectors uniformly distributed in $K$.

In the sequel, if $\mu$ is an isotropic probability measure on $\R^n$,
$f_\theta$ will denote the density of the random variable $\langle
X,\theta\rangle$ with $X$ distributed according to $\mu$. $\gamma$
will denote the density of a standard Gaussian, {\it
  i.e.}, $$\gamma(t) = \frac{1}{\sqrt{2\pi}} e^{-\frac{t^2}{2}}.$$

The following lemma is very well known
\begin{lem}\label{gaussian}
For every $t\geq 1$
$$
\frac{\gamma(t)}{2t}\leq\int_t^\infty\gamma(s) \, ds\leq\frac{2\gamma(t)}{t}.
$$
\end{lem}

The letters $c, c', C, C', c_1,c_2,\ldots$ will denote positive absolute constants, whose value may change from line to line.

Now, let us mention the central limit theorem in the form we will use it to prove Theorem \ref{THM mean
  width} in Section \ref{SEC expected value of mean width}. 

Klartag's central limit theorem for isotropic measures, combined with an
argument by Sodin \cite{key-S} gives the following:

\begin{thrm}[\cite{key-K3}, Theorem 1.4]
  Let $n\geq 1$ be an integer and let $X$ be a random vector in $\R^n$ with an isotropic, log-concave density. Then there exists $\Theta\subseteq S^{n-1}$ with $\sigma_{n-1}(\Theta) \geq 1-Ce^{-\sqrt{n}}$ such that for all 
  $\theta \in \Theta$, the real valued random variable $\skp{X,\theta}$ has a density $f_{\theta}:\R^n \to [0,\infty)$ with the following properties:
  \begin{enumerate}
  \item $\int_{-\infty}^{\infty} \abs{f_{\theta}(t) - \gamma(t)} \, dt \leq \frac{1}{n^{\kappa}}$,
  \item For all $\abs{t}\leq n^{\kappa}$ we have $\abs{\frac{f_{\theta}(t)}{\gamma(t)} - 1} \leq \frac{1}{n^{\kappa}}$.
  \end{enumerate}
  Here, $C,\kappa>0$ are universal constants. 
\end{thrm} 

In the case of symmetric $X$, those results give
$\kappa=\frac{1}{24}$, obtaining that there exists $\Theta\subseteq S^{n-1}$ with $\sigma_{n-1}(\Theta) \geq 1-Ce^{-\sqrt{n}}$ such that for any $\theta\in\Theta$,
  \begin{equation}\label{EQU comparison estimate}
    \abs{\frac{f_{\theta}(t)}{\gamma(t)} - 1} \leq \frac{C'}{n^{\frac{1}{24}}},
  \end{equation}
whenever $\abs{t} < c n^{\frac{1}{24}}$.  

Let us also introduce some notation and results we will need to prove
the estimates for the average of the distribution functions of the 1-dimensional
marginals.

A convex function $M:[0,\infty)\to[0,\infty)$ with $M(0)=0$ and
$M(t)>0$ for $t>0$ is called an Orlicz function (see for instance \cite{key-KR} or
\cite{key-RR}).

Let $X$ be a random vector in $\R^n$. For every $\theta\in S^{n-1}$ we
define an Orlicz function $M_{\theta}$ by
$$
  M_{\theta}(t) = \int_0^s \int_{\{\frac{1}{t}\leq
    \abs{\skp{X,\theta}}\}}\abs{\skp{X,\theta}} \,d\mathbb P \, dt. 
$$ 

This Orlicz function was used in \cite{key-AGP} to study the expected
value of the support function of a random polytope in the direction
$\theta$ since if $M_\theta$ is the Orlicz function associated to a
random vector uniformly distributed on an isotropic body $K$ and $K_N$
is a random polytope on $K$ then we have (see \cite{key-AGP}, Corollary 2.2 )
$$
\E h_{K_N}(\theta)\sim \inf\left\{s>0\,:\,M_\theta\left(\frac{1}{s}\right)\leq\frac{1}{N}\right\}.
$$
The following proposition was obtained:

\begin{prop}[\cite{key-AGP}, Proposition 4.3] \label{PRO david joscha}
 Let $K$ be a symmetric convex body in $\R^n$ of volume 1. Let $s>0$, $\theta\in S^{n-1}$ and $M_{\theta}$ be the Orlicz function associated to the random variable
 $\skp{X,\theta}$, where $X$ is uniformly distributed in $K$. Then,
   \begin{equation}\label{EQU equivalent representation of integral over M}
       \int_{S^{n-1}} M_{\theta}\left(\frac{1}{s}\right) \, d\mu(\theta) = \int_{K} M_{\skp{\theta,e_1}}\left(\frac{\norm{x}_2}{s}\right) \,dx, 
   \end{equation}
 where $M_{\skp{\theta,e_1}}$ is the Orlicz function associated to the random variable $\skp{\theta,e_1}$ with $\theta$ uniformly distributed on $S^{n-1}$.  For any 
 $s\leq \norm{x}_2$
   \begin{equation}\label{EQU formula for M in the general case}
       M_{\skp{\theta,e_1}}\left(\frac{\norm{x}_2}{s}\right) = \frac{2w_{n-1}}{nw_n} \int_0^{\cos^{-1}(\frac{s}{\norm{x}_2})} \frac{\sin^n y}{\cos^2 y} \, dy,
   \end{equation}
 and $0$ otherwise.      
\end{prop} 

With this representation the existence of some directions for which $\mathbb E h_{K_N}(\theta) \geq C \sqrt{\log N} L_K$ holds was established. Using this very same approach we are going to prove Theorem \ref{SupergaussianAverage} and Theorem \ref{SubgaussianAverage} in Section \ref{SuperAndSubGaussianAverage}.

\section{Expected value of the mean width of a random polytope}\label{SEC expected value of mean width}

In this section we are going to prove Theorem \ref{THM mean width}. It is a direct consequence of
the following Theorem, which will fill the gap left by the results already proved in
\cite{key-DGT}.
\begin{thrm}
Let $K$ be a symmetric  isotropic convex body and $K_N$ a random polytope in
$K$, with $n\leq N\leq
n^\delta$. There exist absolute constants $c,C$ and a set
$\Theta\subseteq S^{n-1}$ with $\sigma(\Theta)\geq 1-Ce^{-\sqrt n}$
such that for every
$\theta\in \Theta$
$$
\E h_{K_N}(\theta)\geq \frac{c}{\sqrt\delta}\sqrt{\log N} L_K.
$$
\end{thrm}

\begin{proof}
First of all, notice that for every $\theta\in S^{n-1}$
\begin{eqnarray*}
M_\theta\left(\frac{1}{s}\right)&=&\int_0^{\frac{1}{s}}
\int_{K\cap\{\abs{\skp{x,\theta}}\geq\frac{1}{t}\}} \abs{\skp{x,\theta}} \, dx \, dt  \\
& \geq & \int_{\frac{1}{2s}}^\frac{1}{s}\frac{1}{t}|\{x\in
K\,:\,|\langle x,\theta\rangle| 
\geq \tfrac{1}{t}\}| \, dt \\
&\geq& \left(\frac{1}{s}-\frac{1}{2s}\right)s|\{x\in
K\,:\,|\langle x,\theta\rangle|\geq2s\}| \\
& = & \frac{1}{2}|\{x\in
K\,:\,|\langle x,\theta\rangle|\geq2s\}|.
\end{eqnarray*}
Thus, if $s_0=t_0 L_K$ such that
$M_\theta\left(\frac{1}{s_0}\right)=\frac{1}{N}$ we have that
$$
\frac{2}{N}\geq |\{x\in
K\,:\,|\langle x,\theta\rangle |\geq 2t_0 L_K\}|=\mathbb P
\{|\langle Y,\theta\rangle|\geq 2t_0\},
$$
where $Y$ is a random variable distributed uniformly on
$\frac{K}{L_K}$. Thus, if for some $t$ we have $\mathbb P
\{|\langle Y,\theta\rangle|\geq t\}>\frac{2}{N}$, then, $t_0\geq
\frac{t}{2}$.
From (\ref{EQU comparison estimate}) we have that if $Y$ is  a log-concave
isotropic (covariance matrix equals the identity) random vector in
$\R^n$, then there exists a subset $\Theta\subseteq S^{n-1}$ with
measure greater than $1-Ce^{-\sqrt n}$ such that for any $\theta\in
\Theta$
$$
\left|\frac{f_\theta(t)}{\gamma(t)}-1\right|\leq
\frac{C^\prime}{n^\frac{1}{24}}\textrm{ when } |t|\leq cn^\frac{1}{24}.
$$
Applying this result to the uniform probability measure on
$\frac{K}{L_K}$ we have that there exists $\Theta\subseteq S^{n-1}$ with
measure greater than $1-Ce^{-\sqrt n}$ such that for any $\theta\in
\Theta$ and ant $0\leq t\leq cn^{\frac{1}{24}}$ we have, using Lemma \ref{gaussian}
\begin{eqnarray*}
\mathbb P \{|\langle Y,\theta\rangle|<t\} & \leq & \left(1+\frac{C^\prime}{n^{\frac{1}{24}}}\right)\left(1-2\int_t^\infty\gamma(s) \, ds\right) \\
& \leq & \left(1+\frac{C^\prime}{n^{\frac{1}{24}}}\right)\left(1-\frac{\gamma(t)}{t}\right),
\end{eqnarray*}
and so
$$
\mathbb P \{|\langle Y,\theta\rangle|\geq t\}\geq 1-\left(1+\frac{C^\prime}{n^{\frac{1}{24}}}\right)\left(1-\frac{\gamma(t)}{t}\right).
$$
Taking $t=\alpha\sqrt{\log N}$ we have that $t\leq cn^{\frac{1}{24}}$,
since $N\leq n^\delta$. Thus
\begin{eqnarray*}
\mathbb P \{|\langle Y,\theta\rangle|\geq\alpha\sqrt{\log
  N}\}&\geq&\frac{e^{-\frac{\alpha^2}{2}\log
    N}}{\sqrt{2\pi}\alpha\sqrt{\log
    N}}-\frac{C^\prime}{n^{\frac{1}{24}}}\left(1-\frac{e^{-\frac{\alpha^2}{2}\log
      N}}{\sqrt{2\pi}\alpha\sqrt{\log N}}\right)\cr
&\geq &\frac{1}{\sqrt{2\pi}\alpha N^{\frac{\alpha^2}{2}}\sqrt{\log
    N}}-\frac{C^\prime}{n^{\frac{1}{24}}}\cr
&\geq &\frac{1}{\sqrt{2\pi}\alpha N^{\frac{\alpha^2}{2}}\sqrt{\log
    N}}-\frac{C^\prime}{N^{\frac{1}{24\delta}}}>\frac{2}{N}\cr
\end{eqnarray*}
whenever $N\geq N_0$ if we take $\alpha^2=\frac{1}{24\delta}$. Thus,
for every $\theta\in\Theta$,
$$
\mathbb E h_{K_N}(\theta)\sim s_0>\frac{1}{2\sqrt{24\delta}}\sqrt{\log
  N}L_K.
$$
\end{proof}
\begin{remark}
In a similar way, using the central limit theorem, one can show that
if $n\leq N\leq n^\delta$
with very high probability $h_{K_N}(\theta)\geq \frac{c}{\sqrt
  \delta}\sqrt{\log N} L_K$ for every $\theta\in \Theta$. However, we decided to include this
weaker result instead since it is enough for our purpose and the proof
of Theorems \ref{SupergaussianAverage} and \ref{SubgaussianAverage}
follow the same idea.
\end{remark}

\section{Supergaussian Estimates}\label{SuperAndSubGaussianAverage}

In this section we use the technique introduced in \cite{key-AGP} to
prove Theorem \ref{SupergaussianAverage}, extending the interval in which the average
of the distribution function behaves in supergaussian
way. 
\begin{proof}[Proof of Theorem \ref{SupergaussianAverage}.]
First of all, notice that for any $s>0$ and any isotropic convex body
$K$ we have 
\begin{eqnarray*}
  \int_{S^{n-1}} M_{\theta}\left( \frac{1}{s} \right) \, d\sigma(\theta) & = & \int_{S^{n-1}} \int_0^{\frac{1}{s}} \int_{K\cap\{\abs{\skp{x,\theta}} \geq \frac{1}{t}\}} \abs{\skp{x,\theta}} \, dx \, dt \, d\sigma(\theta) \\
  & \leq & \int_{S^{n-1}} \frac{1}{s} \int_{K\cap\{\abs{\skp{x,\theta}} \geq s\}} \abs{\skp{x,\theta}} \, dx \, d\sigma(\theta) \\
  & \leq & \int_{S^{n-1}} \frac{1}{s} L_K | \{ x \in K : \abs{\skp{x,\theta}} \geq s \}|^{\frac{1}{2}} \, d\sigma(\theta) \\
  & \leq & \frac{L_K}{s} \left( \int_{S^{n-1}} | \{ x \in K : \abs{\skp{x,\theta}} \geq s \}| \, d\sigma(\theta) \right)^{\frac{1}{2}}.
\end{eqnarray*}
Thus, taking $s=tL_K$ we have that for any $t>0$ 
  $$
    \left( \int_{S^{n-1}} |\{ x\in K : \abs{\skp{x,\theta}} \geq tL_K \}| \, d\sigma(\theta) \right)^{\frac{1}{2}} \geq t \int_{S^{n-1}} M_{\theta}\left( \frac{1}{tL_K} \right) \, d\sigma(\theta).
  $$
Using the representation of the average of $M_{\theta}$ as an integral
on $K$ (Proposition \ref{PRO david joscha}) and the lower bound obtained in the proof of Theorem 4.2 in \cite{key-AGP}, we obtain that for
every positive $t$
  \begin{eqnarray*}
    t \int_{S^{n-1}} M_{\theta}\left( \frac{1}{tL_K} \right) \, d\sigma(\theta) & \geq &  t \int_{K\setminus 2tL_KB_2^n} \frac{cw_{n-1}}{n w_n}\frac{\norm{x}_2}{tL_K} e^{-\frac{Cnt^2L_K^2}{\norm{x}_2^2}} \, dx.
  \end{eqnarray*} 
Using the small ball probability estimate from \cite{key-P2}, there
exists an absolute constant $c_1$ such that $|K\backslash c_1\sqrt
nL_K|\geq \frac{1}{2}$. Thus, if $0<t\leq \frac{c_1}{2} \sqrt{n}$   
  \begin{eqnarray*}
    \left( \int_{S^{n-1}} |\{ x\in K :  \abs{\skp{x,\theta}}>tL_K \}| \, d\sigma(\theta) \right)^{\frac{1}{2}} & \geq & 
    \int_{K\setminus c_1\sqrt{n}L_KB_2^n} \frac{cw_{n-1}}{n w_n}\frac{\norm{x}_2}{L_K} e^{-\frac{Cnt^2L_K^2}{\norm{x}_2^2}} \, dx \\
    & \geq & c_2 e^{-Ct^2}\geq e^{-c_3t^2},
  \end{eqnarray*} 
if $1\leq t\leq \frac{c_1}{2}\sqrt n$.
Thus, taking $c=\max\{\frac{2}{c_1},2c_3\}$ we obtain that for every isotropic
convex body $K$ and every $1\leq t\leq \frac{\sqrt{n}}{c}$
  $$
    \int_{S^{n-1}} |\{ x\in K :  \abs{\skp{x,\theta}} > t L_K \}| \, d\sigma(\theta) \geq e^{-ct^2}.
  $$
\end{proof}

As a consequence, we obtain the following estimate of the measure of
the set of directions verifying a supergaussian estimate for a
particular $t$. However, oppositely to what will happen for a
subgaussian estimate, the estimate of the measure of this set of directions will be
very small for big values of $t$.

\begin{cor}
There exist absolute constants $c,c^\prime$ such that for every
isotropic symmetric convex body $K\subseteq\R^n$, and every $1\leq
t\leq \frac{\sqrt n}{c}$ the set of
directions verifying the supergaussian estimate
$$
|\{x\in K\,:\,|\langle x,\theta\rangle|\geq tL_K\}|\geq e^{-c^2t^2}
$$
has measure greater than $ e^{-c^\prime t^2}$.
\end{cor}

\begin{proof}
By Theorem \ref{SupergaussianAverage}
$$
\int_{S^{n-1}}|\{x\in K\,:\,|\langle x,\theta\rangle|<tL_K\}| \, d\sigma(\theta) \leq 1-e^{-c^2t^2}.
$$
Thus, by Markov's inequality, we have
$$
\sigma \left \{\theta\in S^{n-1}\,:\,|\{x\in K\,:\,|\langle x,\theta\rangle|<tL_K\}|>1-e^{-{c^\prime}^2t^2} \right \} \leq\frac{1-e^{-c^2t^2}}{1-e^{-{c^\prime}^2t^2}}.
$$
Consequently,
$$
\sigma\left \{\theta\in S^{n-1}\,:\,|\{x\in K\,:\,|\langle
x,\theta\rangle|\geq tL_K\}|< e^{-{c^\prime}^2t^2} \right \} \leq \frac{1-e^{-c^2t^2}}{1-e^{-{c^\prime}^2t^2}},
$$
and so
\begin{eqnarray*}
\sigma \left \{\theta\in S^{n-1}\,:\,|\{x\in K\,:\,|\langle
x,\theta\rangle|\geq tL_K\}|\geq
e^{-{c^\prime}^2t^2} \right \}&\geq&1-\frac{1-e^{-c^2t^2}}{1-e^{-{c^\prime}^2t^2}}\cr
&=&\frac{e^{-c^2t^2}(1-e^{-({c^\prime}^2-c^2)t^2})}{1-e^{-{c^\prime}^2t^2}}.
\end{eqnarray*}
Choosing ${c^\prime}^2=2c^2$ we obtain the result.
\end{proof}

\section{Subgaussian Estimates}\label{SubgaussianEstimates}

In this section we will extend the interval in which the average of
the distribution function verifies a subgaussian estimate. As a
consequence we will obtain an estimate of the measure of the
directions verifying a subgaussian estimate in some interval.

\begin{proof}[Proof of Theorem \ref{SubgaussianAverage}.]
As we have seen in the proof of Theorem \ref{THM mean width} for every $\theta\in S^{n-1}$, we have
\begin{eqnarray*}
M_\theta\left(\frac{1}{s}\right)&\geq&\frac{1}{2}|\{x\in
K\,:\,|\langle x,\theta\rangle|\geq2s\}|.
\end{eqnarray*}
Thus, using the upper bound shown in \cite{key-AGP}
$$
F(t)\leq
2\int_{S^{n-1}}M_\theta\left(\frac{2}{tL_K}\right) \, d\sigma(\theta)\leq\frac{4\omega_{n-1}}{n\omega_ntL_K} \int_K |x|e^{-\frac{(n-1)t^2L_K^2}{2|x|^2}} \, dx.
$$
By Paouris' concentration of measure result, there exist constants
such that for every $\gamma>c_1$, $|K\backslash \gamma\sqrt nL_K
B_2^n|\leq e^{-c_2\gamma\sqrt n}$. Hence,
\begin{eqnarray*}
\int_K |x|e^{-\frac{(n-1)t^2L_K^2}{2|x|^2}} \, dx
&=& \int_{K\cap\gamma\sqrt nL_K
B_2^n} |x|e^{-\frac{(n-1)t^2L_K^2}{2|x|^2}} \, dx+\int_{K\backslash\gamma\sqrt nL_K
B_2^n} |x|e^{-\frac{(n-1)t^2L_K^2}{2|x|^2}} \, dx\cr
&\leq& \gamma\sqrt n L_K e^{-\frac{ct^2}{\gamma^2}}+(n+1)L_K|K\backslash \gamma\sqrt nL_K
B_2^n|\cr
&\leq&\gamma\sqrt n L_K e^{-\frac{ct^2}{\gamma^2}}+(n+1)L_Ke^{-c_2\gamma\sqrt n}.
\end{eqnarray*}
Consequently,
\begin{eqnarray*}
F(t)&\leq&\frac{C}{t}\left(\gamma e^{-\frac{ct^2}{\gamma^2}}+\sqrt n e^{-c_2\gamma\sqrt
    n}\right)\cr
&\leq&\frac{C}{t}\left(\gamma e^{-\frac{ct^2}{\gamma^2}}+e^{\sqrt
    n}e^{-c_2\gamma \sqrt n}\right)\cr
&\leq&\frac{C}{t}\left(\gamma e^{-\frac{ct^2}{\gamma^2}}+e^{(1-c_2\gamma) t^2}\right)
\end{eqnarray*}
if $1\leq t\leq n^\frac{1}{4}$ and  $\gamma$ is a
constant big enough. Thus, we obtain that for $t$ in the aforementioned interval
$$
F(t)\leq\frac{C}{t}e^{-ct^2}\leq e^{-c't^2}.
$$
Notice that if $R(K)\leq C\sqrt n L_K$, then we do not need to split
the integral as the sum of two integrals and we obtain that
$$
F(t)\leq\frac{C}{t}e^{-ct^2}\leq e^{-c^\prime t^2}
$$
for every $1\leq t\leq C\sqrt n$.
\end{proof}

As a consequence we will obtain that if we consider an interval
$t_0\leq t\leq n^\frac{1}{4}$ with $t_0$ big, then the measure of the
directions that are subgaussian in such interval is big. It is
explicitely stated in the following:

\begin{cor}
Let $K\subseteq \R^n$ be an isotropic symmetric convex body and let
$1\leq t_0\leq n^\frac{1}{4}$. Then the set of directions $\theta\in
S^{n-1}$ that verify
$$
|\{x\in K\,:\,|\langle x,\theta\rangle|\geq tL_K\}|\leq e^{-c^2t^2}
\textrm{ for every }t_0\leq t\leq n^\frac{1}{4}
$$
has measure greater than $1-e^{-{c^\prime}^2 t_0^2}$.\\
Furthermore, if $K$ has small diameter we can take $t_0\leq C\sqrt n$
and $t_0\leq t\leq C\sqrt n$.
\end{cor}

\begin{proof}
Applying Markov's inequality in Theorem \ref{SubgaussianAverage} we
have that for every $1\leq t\leq n^\frac{1}{4}$
$$
\sigma\left \{\theta\in S^{n-1}\,:\,|\{x\in K\,:\,|\langle
x,\theta\rangle|\geq tL_K\right \}\geq e^{-{c^\prime}^2 t^2}\}\leq e^{-(c^2-{c^\prime}^2)t^2}.
$$
Taking ${c^\prime}^2=\frac{1}{2}c^2$ we have that there exists a
constant $c_0$ such that for every $1\leq t\leq n^\frac{1}{4}$
$$
\sigma \left \{\theta\in S^{n-1}\,:\,|\{x\in K\,:\,|\langle
x,\theta\rangle|\geq tL_K \right \}|\geq e^{-{c_0}^2 t^2}\}\leq e^{-c_0^2t^2}.
$$
Take now $t_0< t_1<\dots <t_I= n^{\frac{1}{4}}$. Then,
$$
\sigma \left \{\theta\in S^{n-1}\,:\,\exists 0\leq i\leq I \textrm{ with }|\{x\in K\,:\,|\langle
x,\theta\rangle|\geq t_iL_K\right \}|\geq e^{-{c_0}^2 t_i^2}\}\leq \sum_{i=0}^Ie^{-c_0^2t_i^2}.
$$
Taking $t_i^2=t_0^2+\lambda i$, this probability is bounded by
$e^{-c_0^2t_0^2}\frac{1}{1-e^{-c_0^2\lambda}}$.
If for every $i$ we have that
$$
|\{x\in K\,:\,|\langle
x,\theta\rangle|\geq t_iL_K\}|\leq e^{-{c_0}^2 t_i^2},
$$
then for every $t_i\leq t\leq t_{i+1}$ we have
\begin{eqnarray*}
|\{x\in K\,:\,|\langle
x,\theta\rangle|\geq tL_K\}|&\leq&|\{x\in K\,:\,|\langle
x,\theta\rangle|\geq t_iL_K\}|\leq e^{-{c_0}^2 t_i^2}\cr
&=&e^{-{c_0}^2 t^2}e^{{c_0}^2 (t^2-t_i^2)}\leq e^{-{c_0}^2
  t^2}e^{{c_0}^2 (t_{i+1}^2-t_i^2)}\cr
&=&e^{-{c_0}^2
  t^2}e^{{c_0}^2 \lambda }.
\end{eqnarray*}
Choosing $\lambda$ a constant smaller than 1 we obtain the result.\\
If $K$ has small diameter the same proof works in the interval
$t_0\leq t\leq C\sqrt n$ with $1\leq t_0\leq C\sqrt n$.
\end{proof}

\proof[Acknowledgements]
This work was done while the authors were postdoctoral fellows at the
Department of Mathematical and Statistical Sciences at University of
Alberta. We would like to thank the department for providing such good
environment and working conditions. We would also like to thank
Prof. Jesus Bastero and Dr. Peter Pivovarov for helpful comments and for reading a preliminary
version of this paper.

\end{document}